\documentclass{blamvart}
\newtheorem{theorem}{Th\'eor\`eme}[section]
\newtheorem{lemma}[theorem]{Lemme}

\theoremstyle{definition}

\theoremstyle{remark}
\newtheorem{remark}[theorem]{Remarque}
\numberwithin{equation}{section}

\title[Sur les points fixes r\'epulsifs]{Sur les points fixes et les
cycles r\'epulsifs au voisinage d'une 
singularit\'e essentielle isol\'ee \`a l'instar de la m\'ethode de renormalisation de Zalcman}

\author{Claudi Meneghin} 
\address{Claudi Meneghin - Fermo Posta Chiasso 1 - CH-6830 Chiasso - Switzerland}
\email{claudi.meneghin@gmail.com}

\begin{document}
\begin{abstract}
Soit $g$ une fonction holomorphe au voisinage 
d'une singularit\'e essentielle isol\'ee $v$: nous prouvons que,
si $g$ y omet une valeur complexe, alors $v$ 
peut \^etre approch\'e par une suite de points fixes
r\'epulsifs de $g$, dont les multiplicateurs divergent \`a $\infty$.
Nous montrons aussi que, si $v$ n'est pas une valeur exceptionnelle au sens de Picard pour $g$, 
alors $v$ peut \^etre approch\'ee par une suite de points 
p\'eriodiques d'ordre deux de $g$, ces cycles \'etant 
r\'epulsifs (avec multiplicateurs divergeant \`a $\infty$) 
si $v$ n'est pas une valeur compl\`etement ramifi\'ee.
\end{abstract}
\maketitle

\font\sdopp=msbm10 
\font\sdoppet=msbm8
\font\cir=wncyb10
\def\ER {\sdopp {\hbox{R}}}
\def\CI {\sdopp {\hbox{C}}}
\def\Cit {\sdoppet {\hbox{C}}}
\bibliographystyle{plain}

\font\cir=wncyb10
\def\iu{{\cir\hbox{Z}}}
\def\ze{{\cir\hbox{Z}}}
\def\pe{{\cir\hbox{P}}}
\def\ef{{\cir\hbox{F}}}

\parindent=8pt
\font\cir=wncyb10
\def\Iu{\cir\hbox{YU}}
\def\Ze{\cir\hbox{Z}}
\def\pe{\cir\hbox{P}}
\def\Ef{\cir\hbox{F}}
\def\CIRC{\mathop{\tt o}\limits}
\def\quan{\vrule height6pt width6pt depth0pt}
\def\QUAN{\ \quan}
\def\BETA{\mathop{\beta}\limits}
\def\GAMMA{\mathop{\gamma}\limits}
\def\VI{\mathop{v}\limits}
\def\UI{\mathop{u}\limits}
\def\VII{\mathop{V}\limits}
\def\WI{\mathop{w}\limits}
\def\ZETA{\mathop{Z}\limits}
\def\W{{{\cal W}\, }^{\circ}}

\font\sdopp=msbm10
\def\ER {\sdopp {\hbox{R}}}
\def\QU {\sdopp {\hbox{Q}}}
\def\CI {\sdopp {\hbox{C}}}
\def\DI {\sdopp {\hbox{D}}}
\def\EN{\sdopp {\hbox{N}}}
\def\ZETA{\sdopp {\hbox{Z}}}
\def\ES{\sdopp {\hbox{S}}}
\def\PI {\sdopp {\hbox{P}}}
\def\M{\hbox{\tt\large M}}
\def\N{\hbox{\boldmath{}$N$\unboldmath}} 
\def\P{\hbox{\boldmath{}$P$\unboldmath}} 
\def\tr{\hbox{\boldmath{}$tr$\unboldmath}} 
\def\f{\hbox{\large\tt f}} 
\def\F{\hbox{\boldmath{}$F$\unboldmath}} 
\def\G{\hbox{\boldmath{}$G$\unboldmath}} 
\def\L{\hbox{\boldmath{}$L$\unboldmath}} 
\def\h{\hbox{\boldmath{}$h$\unboldmath}} 
\def\e{\hbox{\boldmath{}$e$\unboldmath}} 
\def\g{\hbox{\boldmath{}$g$\unboldmath}} 
\def\u{\hbox{\boldmath{}$u$\unboldmath}} 
\def\v{\hbox{\boldmath{}$v$\unboldmath}} 
\def\U{\hbox{\boldmath{}$U$\unboldmath}} 
\def\V{\hbox{\boldmath{}$V$\unboldmath}} 
\def\id{\hbox{\boldmath{}$id$\unboldmath}} 
\def\alph{\hbox{\boldmath{}$\alpha$\unboldmath}} 
\def\bet{\hbox{\boldmath{}$\beta$\unboldmath}} 
\def\gam{\hbox{\boldmath{}$\gamma$\unboldmath}} 
\def\pphi{\hbox{\boldmath{}$\varphi$\unboldmath}} 
\def\ppsi{\hbox{\boldmath{}$\psi$\unboldmath}} 
\def\Ppsi{\hbox{\boldmath{}$\Psi$\unboldmath}} 
\def\brevve{}
\def\labelle #1{\label{#1}}
\def\quadras #1{ \hbox{\bf[}{#1}
\hbox{\bf]}}
\def\cal{}
\font\sdopp=msbm10 scaled \magstep1
\font\sdoppon=msbm6 scaled \magstep1
\def\CI {\sdopp {\hbox{C}}}
\def\CIP {\sdoppon {\hbox{C}}}

\section{Introduction} 

Soit $f$ une fonction analytique dans une r\'egion du plan complexe. 
Soit $f^{\circ p}$ l'it\'er\'e p-i\`eme de $f$, c'est-\`a-dire, 
d\'efinissons $f^{\circ 2} (z) := f(f(z))$ et, pour $p\in\EN$, $f^{\circ p} (z) := f(f^{\circ p-1}(z))$
(lorsque ces quantit\'es sont bien d\'efinies).

Un point $z_0$ est dit un {\it point p\'eriodique de p\'eriode} $p$ de $f$ si $p$ est un entier naturel 
tel que $f^{\circ p} (z_0)$ est bien d\'efini et $f^{\circ p} (z_0) = z_0$. 

Si p est le plus petit nombre naturel tel que $f^{\circ p} (z_0) = z_0$, alors le point $z_0$ est appel\'e un 
{\it point p\'eriodique de p\'eriode minimale} $p$. 
L'ensemble $$\{z_0, f(z_0), f^{\circ 2}(z_0), \dots, f^{\circ (p-1)}(z_0)\}$$ 
est appel\'e un {\it cycle de points p\'eriodiques}. 
La valeur $A := (f^{\circ p})'(z_0)$ est appel\'ee le {\it multiplicateur} du point p\'eriodique $z_0$.
Un cycle de p\'eriode minimale $1$ est appelé un {\it point fixe}. Un cycle (y compris les points fixes) est dit 
{\it r\'epulsif} si $| A |> 1$.

Bergweiler a prouv\'e \cite{bergweiler_int} que, pour tout $n>1$, toute fonction enti\`ere transcendante a une suite de $n-$cycles r\'epulsifs \`a multiplicateurs divergeant \`a $\infty$. Dans cette note, nous prouvons un r\'esultat en corr\'elation:
\begin{quote}
\it
Soient $ v\in\CI $, ${\cal W}$ 
un voisinage ferm\'e de $v$, $g$ une fonction holomorphe
sur ${\cal W}\setminus \{v\}$, ayant une
singularit\'e essentielle \`a $v$. 
S'il existe une valeur complexe $\alpha $ omise par $g$ 
sur ${\cal W}\setminus \{v\}$, alors il existe une suite $\{q_n\}\to v$ de points fixes r\'epulsifs de $g$, dont les multiplicateurs divergent \`a $\infty$ (th\'eor\`eme \ref{principal}).
En outre (th\'eor\`eme \ref{principal2}), si $v$ n'est pas une valeur exceptionnelle au sens de Picard pour $g$, 
alors $v$ peut \^etre approch\'ee par une suite de points 
p\'eriodiques d'ordre deux de $g$, ces cycles \'etant 
r\'epulsifs (avec multiplicateurs divergeant \`a $\infty$) 
si $v$ n'est pas une valeur compl\`etement ramifi\'ee.
\end{quote}
La preuve de ce r\'esultat utilise de fa\c con essentielle la m\'ethode de renormalisation de Zalcman (voir e.g.$\!$ \cite{bzalcman}):
cette technique a d\'esormais une vaste gamme d'applications; par exemple Schwick \cite{schwick} l'a utilis\'ee par simplifier la preuve de Baker \cite{baker_rep} que l'ensemble de Julia d'une fonction enti\`ere transcendante est la fermeture de l'ensemble des 
points p\'eriodiques r\'epulsifs. Pourtant, dans cet article, nous ne l'utiliserons pas par enqu\^eter sur l'ensemble de Julia des
fonctions analytiques, mais sur leur dynamique au voisinage d'une singularit\'e essentielle isol\'ee.
Notons en outre que les th\'eor\`emes \ref{principal} et \ref{principal2} sont de toute fa\c con purement locaux: 
en d'autres mots ils ne d\'emandent que les fonctions concern\'ees soyent analytiques sur $\CI$. 

L'adaptation de la m\'ethode de Zalcman utilise un th\'eor\`eme de Lehto et Virtanen sur la croissance de la d\'eriv\'ee sph\'erique 
au voisinage d'une singularit\'e essentielle isol\'ee (th\'eor\`eme {\ref{rapp}}, voir \cite{lehto&virtanen, lehto}) ainsi que une
g\'en\'eralisation (lemme {\ref{metric}}) du {\it lemme de l'espace m\'etrique} de Gromov (voir \cite{gromov}, p.256).

La preuve du th\'eor\`eme \ref{principal} consiste essentiellement \`a composer $g$ \`a la source de avec une suite de contractions bien choisies, ce qui 
permet de construire une fonction holomorphe enti\`ere limite; on conclut en appliquant les
th\'eor\`emes de Picard, Hurwitz (voir \cite {berteloot}, p.8) et des quatre 
valeurs compl\`etement ramifi\'ees (voir \cite{bergdyn}, th.\,29 et 30).


\section{Pr\'eliminaires}

Rappelons d'abord qu'une fonction m\'eromorphe $f$ est dite {\it faiblement normale} sur un domaine ${\cal D}$ si, pour chaque sous-domaine simplement 
connexe $G\subset {\cal D} $, la famille $\{f\circ S\}_{S\in \hbox{\tt\small Aut}(G)}$, index\'ee sur les automorphismes de $G$, est normale; soit 
$f^{\sharp}$ la d\'eriv\'ee sph\'erique de $f$. On a:
\begin{theorem}
Soit $v\in\CI$, ${\cal W}$ un voisinage de
$v$ in $\CI$; soit $g$ une fonction holomorphe 
${\cal W}\setminus\{v\}\to\CI$, ayant une singularit\'e essentielle isol\'ee \`a $v$ et omettant une valeur complexe $\alpha$ au voisinage de $v$.
Alors
$ \limsup_{z\to v} \vert z-v\vert\cdot g^{\sharp}(z)=\infty $.
\labelle{rapp}
\end{theorem}
{\bf D\'emonstration:} gr\`ace au th\'eor\`eme d'Iversen, la valeur $\alpha$ est une valeur asymptotique de $g$; 
gr\^ace au th\'eor\`eme 2 en \cite{lehto&virtanen} (voir aussi 
\cite{lehto&virtanen2}, point 7), $g$ n'est pas faiblement normal en ${\cal W}\setminus\{v\}$; la th\`ese s'ensuit alors du 
th\'eor\`eme 2 de \cite{lehto&virtanen2}.
\QUAN
\vskip0.2cm
Le lemme suivant est une version renforc\'ee du
{\sf lemme de l'espace m\'etrique} (voir \cite{gromov}, p. 256), n\'ecessaire pour la d\'emonstration du lemme de renormalisation
\ref{zalcman}:

\begin{lemma} 
Soit $(X,d)$  un espace m\'etrique  complet,
$Y\subset X$ un sous-ensemble de $X$ 
tel que $X\setminus \overline Y\not=\emptyset$
et 
$M:X\rightarrow \ER^+$ 
une fonction localement born\'ee sur $X\setminus \overline Y$.
Soit $\sigma>0$:  alors, 
pour tout $u\in X$ 
tel que 
$d(Y, u)> {2}/{\sigma M(u)}$
il existe
$w\in X$
tel que:
\begin{eqnarray}
&\ & d(u,w) \leq \quadras{\sigma M(u)}^{-1}
\labelle{metric1}\\
&\ & M(w) \geq  M(u)
\labelle{metric2}\\
&\ &
\overline{D}\left(w, \quadras{\sigma M(w)}^{-1}\right)
\cap Y  =  \emptyset 
\labelle{metric3}\\
&\ &  d(x,w)\leq \quadras{\sigma M(w)}^{-1}
 \Rightarrow 
M(x)\leq 2 M(w)
\labelle{metric4}.
\end{eqnarray}
\labelle{metric}
\end{lemma}
{\bf D\'emonstration:} 
$ u\in X\setminus \displaystyle\overline{Y}$ \'etant donn\'e,
tel que $d(Y, u)> {2}/{\sigma M(u)}$, supposons par l'absurde qu'il n'existe pas le $w$ du lemme. 
Alors $v_0:= u$ n'est pas convenable en tant que choix de $w$. 
Par contre, $w=u$ v\'erifie automatiquement 
(\ref{metric1}), (\ref{metric2}) et (\ref{metric3}), donc il
doit violer la condition (\ref{metric4}). Ainsi on peut trouver $v_1\in X$
tel que $M(v_1)>2 M(v_0)$ et  $d(v_1,v_0)\leq\quadras{\sigma M(v_0)}^{-1}$.

Par cons\'equent, $v_1$ aussi v\'erifie l'hypoth\`ese du lemme, car 
$
d(Y, v_1)
\geq
d(Y, v_0)-d(v_1, v_0)
\geq  2\quadras{\sigma M(v_0)}^{-1} 
- \quadras{\sigma M(v_0)}^{-1}=
\quadras{\sigma M(v_0)}^{-1}
>
2\quadras{\sigma M(v_1)}^{-1} 
$.

Par induction, 
on peut ainsi construire une suite 
$\{v_n\}$ telle que $v_0=u$, 
$ M(v_{n+1})>2 M(v_n)$ mais 
$d(v_n,v_{n+1})\leq \quadras{\sigma M(v_n)}^{-1}$; par cons\'equent, 
\begin{eqnarray}
M(v_{n})&\geq&
2^{n}M(v_0)
\labelle{croissance}\\
d(v_n, v_{n+1})&\leq& \quadras{2^n\sigma M(v_0)}^{-1}.
\labelle{borne}
\end{eqnarray}
Ainsi, 
\begin{eqnarray*}
d(v_{n},v_{n+k})
\leq
\sum_{l=n}^{n+k-1}
d(v_{l},v_{l+1})
<\sum_{l=n}^{\infty}
d(v_{l},v_{l+1})
\leq\\
\leq
\left[\sigma M(v_0)\right]^{-1}
\sum_{l=n}^{\infty} {2} ^{-l}
=
\left[\sigma M(v_0)\right]^{-1}
{2}^{1-n}
,\end{eqnarray*} 
ce qui entra\i ine que la suite $\{v_n\}$ est de Cauchy. 

En soit $\lambda$ la valeur limite:
en rappelant que $v_0=u$, on a 
$$d(u, \lambda)=\lim_{k\to\infty}d(v_0, v_k)
\leq 2\quadras{\sigma M(u)}^{-1},$$  
ce qui entra\^\i ne 
$d(Y,\lambda)\geq d(Y,u)-d(u,\lambda )>0$, 
car $d(Y,u)> 2\quadras{\sigma M(u)}^{-1}$. 
Ainsi $\lambda\not\in\overline{Y}$; par contre, gr\^ace \`a (\ref{croissance}),
$M$ n'est pas born\'e
au voisinage de $\lambda$: c'est une contradiction.
\QUAN
\vskip.5cm
Le lemme suivant 'renormalise', \`a l'instar de la m\'ethode de 
Zalcman, une fonction holomorphe 
$g$ au voisinage d'une singularit\'e essentielle isol\'ee. 
En effet, on proc\'edera en composant $g$ avec un suite de 
contractions; pourtant, on ne sera pas concern\'e avec une famille 
pas normale de fonctions m\'eromophes, mais avec une seule fonction \`a singularit\'e essentielle isol\'ee.

\begin{lemma}
Soient $v\in\CI$, ${\cal W}$ un voisinage ferm\'e
de $v$, $g$ une fonction holomorphe
sur ${\cal W}\setminus \{v\}$, ayant une
singularit\'e essentielle \`a $v$.
Alors il existe des
suites $\{v_n\}\to v$, $\{r_n\}\subset\ER^+$, 
avec $\{r_n\}\to 0$, telles que
$\{g (v_n+r_n z)\}$ converge uniform\'ement sur 
tout compact de $\CI$ vers une application holomorphe non 
constante $h:\CI\to\CI$.
En outre, si $g$ ne prend pas la valeur $\alpha\in\CI$ sur
${\cal W}\setminus \{v\}$, alors $h$ ne prend
pas la valeur $\alpha$ non plus. 
\labelle{zalcman}
\end{lemma}
{\bf D\'emonstration:} gr\^ace au th\'eor\`eme \ref{rapp},
on peut trouver des suites $\{\lambda_n\}\rightarrow +\infty$ 
en $\ER$ et $\{\xi_n\}\rightarrow v$ 
en ${\cal W}\setminus\{v\}$ telles que 
\begin{equation}
\vert  
\xi_n - v \vert\cdot g^{\sharp}(\xi_n)
=
\lambda_n
.
\labelle{lambda}
\end{equation}
Pour tout $n$, le lemme \ref{metric} est applicable 
\`a ${\cal W}$ avec la m\'etrique euclid\'eenne, 
$Y=\{v\}$,
$M(x)=
g^{\sharp} (x)$, 
$u=\xi_n$
et $\sigma=3/\lambda_n$: en effet, gr\^ace \`a (\ref{lambda}), 

$$
d(Y, u)=
\vert  
\xi_n - v \vert =
\frac{\lambda_n}
{g^{\sharp}(\xi_n)}
=
\frac{3}
{\sigma g^{\sharp}(\xi_n)}
>
\frac{2}
{\sigma g^{\sharp}(\xi_n)}
=
\frac{2}
{\sigma M(u)}
\, .
$$

On obtient $v_n\in {\cal W}$ 
tel que:
\begin{eqnarray}
&\ &
\vert \xi_n - v_n \vert
\leq \frac{\lambda_n}{3g^{\sharp}(\xi_n)}
=
\vert \xi_n - v\vert
\labelle{tt1}
\\
&\ &
g^{\sharp}(v_n)     
\geq 
g^{\sharp}(\xi_n)
\labelle{tt2}
\\
&\ &
\overline{\DI}(v_n, \frac{\lambda_n}{3g^{\sharp}(v_n)})\cap\{v\}=\emptyset 
\labelle{tt3}
\\
&\ &
\vert x-v_n     \vert
\leq 
\displaystyle
\frac{\lambda_n}
{3 g^{\sharp} (v_n)}
\Rightarrow
g^{\sharp}(x)\leq 2 g^{\sharp}(v_n)\, .
\labelle{tt4}
\end{eqnarray}

Par cons\'equent, $\, v_n\to v\, $ car, gr\^ace \`a (\ref{tt1}), 
$
\vert v_n - v\vert
\leq 
\vert \xi_n - v_n\vert+
\vert \xi_n - v\vert
\leq
\frac{4}{3}\,\vert \xi_n - v\vert
$.

Posons maintenant
$r_n:=\quadras{3g^{\sharp} (v_n)}^{-1}$ et
$h_n(z):=g(v_n +r_n z)$:
on voit sur (\ref{tt3}) que 
\begin{eqnarray*}
z\in \DI(0,\lambda_n)&\Rightarrow& 
v_n +r_n z\in
\DI\left(v_n, \frac{\lambda_n}{\quadras{3g^{\sharp} (v_n)}}\right)
\\
&\subset& 
\W\setminus\{v\},
\end{eqnarray*}
pour $n$ suffisamment grand (ce que nous
sous-entendrons dans la suite).

Gr\^ace \`a (\ref{tt2}) et (\ref{tt3}),
chaque $h_n$ est bien d\'efini sur
$ \DI(0,\lambda_n)$.

La famille $\{h_n\}$ est normale,  car, 
gr\^ace \`a (\ref{tt4}): 
\begin{equation}
z\in \DI(0,\lambda_n)\Rightarrow h_n^{\sharp}(z)=
g_n^{\sharp}\left(v_n+\frac{z}
{\quadras{3g^{\sharp} (v_n)}}\right)
\,\quadras{3g^{\sharp} (v_n)}^{-1}
\leq 2.
\label{douze}
\end{equation} 

Gr\^ace au th\'eor\`eme d'Ascoli, et au fait que 
$\{\lambda_n\}\to\infty$, on peut extraire de $\{h_n\}$
une sous-suite uniform\'ement convergente 
(que nous appellerons une fois de plus $\{h_n\}$) sur tout compact de
$\CI$ vers une application m\'eromorphe limite 
$h$.
Cette application jouit de la propri\'et\'e que
$${h}^{\sharp}(0) 
=\lim_{n\to\infty}
  h_n^{\sharp}(0) 
=1/3,$$ 
ce qui prouve qu'elle n'est pas constante.
Comme, pour tout $n$, $h_n$ ne prend
pas la valeur $\infty$ sur $\DI(0,\lambda_n)$, 
gr\^ace au lemme de Hurwitz, $h$ ne prend pas la valeur $\infty$: c'est donc 
holomorphe.

Finalament, si $g$ ne prend pas la valeur $\alpha$ sur
${\cal W}\setminus \{v\}$, alors, pour tout $n$, $h_n$ ne prend
pas la valeur $\alpha$ sur $\DI(0,\lambda_n)$ et, gr\^ace une fois de plus au lemme de Hurwitz, 
$h$ ne prend pas la valeur $\alpha$.
\QUAN

\begin{remark}

On voit facilement sur (\ref{douze}) que pour tout $z\in\CI$, 
${h}^{\sharp}(z) \leq 2$. Nous n'utiliserons pas cette propri\'et\'e dans cet article.
\end{remark}

Enfin, on aura besoin d'une extension immediate du th\'eor\`eme 31 de \cite{bergdyn}:

\begin{theorem}
Pour tout $\zeta\in\CI$, toute application holomorphe non constante de 
$h:\CI\setminus\{\zeta\}\to \CI\setminus\{\zeta\}$ n'a pas 
de valeurs compl\`etement ramifi\'ees.
\labelle{punctured}
\end{theorem}
{\bf D\'emonstration:} 
Si $h$ est transcendante, c'est une cons\'equence immediate des deux th\'eor\`emes fondamentaux 
de la th\'eorie de Nevanlinna, voir le th\'eor\`eme 31 de \cite{bergdyn}.
Si par contre $h$ est un polyn\^ome, alors on a $h(z)=\alpha(z-\zeta)+\zeta$ car sinon $h$
prendrait la valeur $\zeta$ pour $z\ne\zeta$, ainsi $h'(z)\equiv\alpha\ne 0$.
\QUAN

\section{Les r\'esultats principaux}
\begin{theorem}
Soient $ v\in\CI $, ${\cal W}$ 
un voisinage ferm\'e de $v$, $g$ une fonction holomorphe
sur ${\cal W}\setminus \{v\}$, ayant une
singularit\'e essentielle \`a $v$. 
S'il existe une valeur complexe $\alpha $ omise par $g$ 
sur ${\cal W}\setminus \{v\}$, alors il existe une suite $\{q_n\}\to v$ de points fixes r\'epulsifs de $g$, dont les multiplicateurs divergent \`a $\infty$.
\labelle{principal}
\end{theorem}
{\bf D\'emonstration:}
\vskip0cm 
{\tt A)} Consid\'erons d'abord le cas $\alpha\not=v$:
gr\`ace au lemme {\ref{zalcman}}, on trouve des suites 
$\{v_{n}\}\to v$ et
$\{r_{n}\}\downarrow 0$ 
telles que
$h_n(z):=\{g(v_{n}+r_{n}z)\}$ converge uniform\'ement sur tout compact de $\CI$
vers une fonction holomorphe enti\`ere non constante $h$ (\`a valeurs en $\CI$).
Gr\^ace au lemme {\ref{zalcman}} $h$ omet la valeur 
$\alpha\not=v$, donc, gr\^ace au th\'eor\`eme de Picard, il prend la valeur $v$, cette valeur n'\'etant pas compl\`etement 
ramifi\'ee, gr\^ace au 
th\'eor\`eme \ref{punctured}.
Ainsi, il existe $z_0\in\CI$ tel que 
\begin{equation}
\begin{cases}
h(z_0)=v\cr 
h^{\prime}(z_0)\not= 0.
\end{cases}
\labelle{branche1}
\end{equation}

Or, $z\mapsto g(v_{n}+r_{n}z)-(v_{n}+r_{n}z)$
converge, apres eventuelle extraction,
vers $h-v$, et $h(z_0)-v=0$,
donc le lemme de Hurwitz 
nous passe une suite de points
$\{z_{n}\}\to z_0$ 
tels que 
$g(v_{n}+r_{n}z_{n} )=(v_{n}+r_{n}z_{n} )$: 
ainsi les points 
$q_{n}:=v_{n}+r_{n}z_{n} $ 
forment une suite $ \{q_{n}\}\to v  $
de points fixes de $g$.

Ces points sont r\'epulsifs (pour $n$ assez grand) car on a, 
d'un c\^ot\'e, gr\^ace au choix de $z_0$ en (\ref{branche1}):
$$(g\circ h)^{\prime}(z_0)
=
h^{\prime}(z_0)
\cdot
g^{\prime}\left(h(z_0)\right)
\not=0  
$$ et de l'autre c\^ot\'e, 
$r_{n}\to 0$, 
$$
r_{n}\cdot
g^{\prime}
(v_{n}+r_{n}z_{n} )
=h_n^{\prime}(z_n)
\to
h^{\prime}(z_0)
,$$
ce qui prouve que, pour $n$ assez grand, les $q_{n}$
sont r\'epulsifs et 
$g^{\prime}(q_{n})\to\infty$.
\vskip.2cm
{\tt B)} Soit maintenant $\alpha=v$, c'est-\`a-dire $g$ omet la 
m\^eme valeur complexe $v$ au voisinage de la singularit\'e essentielle $v$.

Appliquons le lemme {\ref{zalcman}} \`a la fonction 
$({g(z)-v})/({z-v})$; c'est correct, 
car cette fonction est holomorphe (\`a valeurs en $\CI$) 
sur ${\cal W}\setminus\{v\}$, a une singularit\'e 
essentielle isol\'ee en $z=v$ et ne prend pas la valeur $0$ sur
${\cal W}\setminus \{v\}$.

On trouve donc $\{v_n\}\to v$ et $\{r_n\}\to 0$ tels que
\begin{equation}
h_n(z):= \frac{g(v_{n}+r_{n}z_{n})-v }{v_{n}+r_{n}z_{n}-v}\longrightarrow 
h(z),
\labelle{pasdepoles1}
\end{equation}
uniform\'ement sur tout compact de $\CI$, o\`u $h$ est une fonction
holomorphe enti\`ere non constante. 
Gr\^ace au lemme de Hrwitz $h$ n'a pas de p\^oles; naturellement $h$ n'a pas de singularit\'es essentielles en $\CI$
donc on a: 
\begin{equation}
\left\vert
\frac{v-v_{n}}{r_{n}}
\right\vert
\longrightarrow 
\infty
\labelle{pasdepoles2}
\end{equation}
en (\ref{pasdepoles1}), car sinon 
$\lim_{n\to\infty}({v-v_{n}})/{r_{n}}$ serait un p\^ole ou une singularit\'e essentielle pour $h$.

Gr\^ace au lemme {\ref{zalcman}}, $h$ ne prend pas la valeur $0$:
gr\^ace au th\'eor\`eme de Picard, $h$ prend alors la valeur $1$. 
Cette valeur n'est pas compl\`etement ramifi\'ee, gr\^ace 
au th\'eor\`eme \ref{punctured}. 
Il existe alors 
$z_0\in\CI$ tel que 
\begin{equation}
\begin{cases}
h(z_0)=1\cr 
h^{\prime}(z_0)\not= 0.
\end{cases}
\labelle{branche10}
\end{equation}
Le lemme de Hurwitz 
nous donne une suite de points
$\{z_{n}\}\to z_0$ 
tels que 
$$\frac{g(v_{n}+r_{n}z_{n} )-v}{(v_{n}+r_{n}z_{n} )-v}=1;$$ 
ainsi les points 
$q_{n}:=v_{n}+r_{n}z_{n} $ 
forment une suite $ \{q_{n}\}\to v  $
de points fixes de $g$. 
Prouvons que ces points sont r\'epulsifs (pour $n$ assez grand). 
On a, par d\'efinition:
$$g(v_{n}+r_{n}z_{n})
=
v+
\left(v_{n}+r_{n}z_{n} - v\right)
\cdot
h_{n}(z).  
$$ 
En d\'erivant en $z=z_n$ et en divisant pour $r_n$, on obtient
\begin{eqnarray*}
g^{\prime}(v_{n}+r_{n}z_{n} )
&=&
h_{n}(z_{n})
+
\left(\frac{v-v_{n}}{r_{n}}+z_n\right)
\,h_{n}^{\prime}(z_{n})\\
\end{eqnarray*}  

Gr\^ace \`a (\ref{pasdepoles2}), $\left\vert
({v-v_{n}})/{r_{n}}
\right\vert
\longrightarrow 
\infty$; comme on a aussi $h_{n}(z_{n})\to h(z_0)$ et 
$h_{n}^{\prime}(z_{n})\to h^{\prime}(z_0)\not=0$ 
(voir (\ref{branche10})), on en tire que 
$g^{\prime}(v_{n}+r_{n}z_{n} )\to\infty$,
ce qui conclut la d\'emonstration.
\QUAN
\vskip0.2cm
Enfin, on a aussi:
\begin{theorem}\labelle{principal2}
Soient $ v\in\CI $, ${\cal W}$ 
un voisinage ferm\'e de $v$, $g$ une fonction holomorphe
sur ${\cal W}\setminus \{v\}$, ayant une
singularit\'e essentielle \`a $v$. 
Si $v$ n'est pas une valeur exceptionelle au sens de Picard de $g$ 
au voisinage de $v$, il existe une suite 
$\{q_n\}\to v$ de 2-cycles de $g$.
En outre, si $v$ n'est pas une valeur compl\`etement ramifi\'ee de
$g$, les $\{q_n\}$ peuvent \^etre choisis r\'epulsifs, avec multiplicateurs divergeant \`a $\infty$. 
\end{theorem}
{\bf D\'emonstration:} 
gr\`ace au lemme {\ref{zalcman}}, on trouve des suites 
$\{v_{n}\}\to v$ et
$\{r_{n}\}\downarrow 0$ 
telles que
$\{h_n(z)\}:=\{g(v_{n}+r_{n}z)\}$ converge uniform\'ement sur tout compact de $\CI$
vers une fonction holomorphe enti\`ere non constante $h$.
Gr\^ace \`a l'hypoth\`ese sur le caract\`ere non 
exceptionnel de la valeur $v$, l'ensemble 
$${\cal S}:= g^{-1}(v)\cap\W\setminus \{v\} $$ 
est infini.
Gr\^ace aux th\'eor\`eme de Picard (ou de fa\c con banale, si $h$ est un polyn\^ome) il existe 
\begin{equation}
\begin{cases}
z_0\in\CI\cr 
w\in{\cal S}
\end{cases}
\labelle{branche0}
\end{equation}
tels que $h(z_0)=w$.
Par continuit\'e, il existe aussi un voisinage ${\cal U}$
de $z_0$ tel que $ h({\cal U})\subset\W\setminus \{v\} $
et, par cons\'equent, les 
$ g^{\circ 2}(v_n+ r_n z)$ sont bien d\'efinis, pour $n$ assez grand, sur ${\cal U}$.
Or,
$z\mapsto\{g^{\circ 2}(v_{n}+ r_{n}z)-(v_{n}+r_{n}z)\}$
converge, apr\`es eventuelle extraction, vers
$g\circ h-v$ uniform\'ement sur tout compact
de ${\cal U}$. 
Comme 
$g\circ h(z_0)-v=g(w)-v=0$,
le lemme de Hurwitz nous passe une suite de points
$\{z_{n}\}\to z_0$ 
tels que 
$g^{\circ 2}(v_{n}+r_{n}z_{n})-(v_{n}+r_{n}z_{n} )=0$,
ainsi les points 
$q_{n}:=v_{n}+ r_{n}z_{n} $ 
forment une suite $ \{q_{n}\}\to v $
de 2-cycles de $g$.
Si, de plus, $v$ n'est pas une valeur compl\`etement ramifi\'ee de
$g$, on peut choisir les $z_0\in\CI$ et $w\in{\cal S}$ en (\ref{branche0}) de fa\c con que
\begin{equation}
\begin{cases}
h(z_0)=w\cr 
h^{\prime}(z_0)\not= 0\cr
g^{\prime}(w)\not= 0.
\end{cases}
\labelle{branche}
\end{equation}

En effet, comme la valeur $v$ n'est pas compl\`etement 
ramifi\'ee, l'ensemble 
$${\cal T}:=\{w\in{\cal S}: g^{\prime}(w)\not= 0\}$$
(et, par cons\'equent, $h^{-1}({\cal T})$) est infini.
L'ensemble 
$$\{z\in h^{-1}({\cal T}): h^{\prime}(z)\not= 0\}$$ 
est de m\^eme infini: cela d\'ecoule, si $h$ est transcendant,
du th\'eor\`eme des quatre valeurs compl\`etement ramifi\'ees; 
si $h$ est un polyn\^ome, alors ceci s'ensuit tout simplement du 
fait que l'ensemble $\{z\vert h^{\prime}(z)= 0\}$ est fini. 
Cela prouve (\ref{branche}): avec ce choix 
de $w$ et $z_0$, on a:
$$(g\circ h)^{\prime}(z_0)
=
h^{\prime}(z_0)
\cdot
g^{\prime}\left(h(z_0)\right)
\not=0 ; 
$$
d'autre c\^ot\'e, 
$r_{n}\to 0$ et
$$
r_{n}\cdot
(g^{\circ 2})^{\prime}
(v_{n}+r_{n}z_{n} )
=(g\circ h_n)^{\prime}(z_n)
\to
(g\circ h)^{\prime}(z_0)
,$$
ainsi, pour $n$ assez grand, les $q_{n}$
sont r\'epulsifs; on a aussi 
$(g^{\circ 2})^{\prime}(q_{n})\to\infty$.
\QUAN

\end{document}